\documentclass{amsart}

\usepackage{amssymb} 
\usepackage{amsmath} 
\usepackage{amsthm} 
\usepackage{graphicx} 
\usepackage{amscd} 
\usepackage{setspace} 

\newtheorem{theorem}{Theorem}[section]

\newtheorem{proposition}[theorem]{Proposition}
\newtheorem{question}[theorem]{Question}

\theoremstyle{remark} 
\newtheorem{example}[theorem]{Example}
\newtheorem{remark}[theorem]{Remark}
\newtheorem{acknowledgements}[theorem]{Acknowledgements}

\theoremstyle{definition} 
\newtheorem{definition}[theorem]{Definition}

\newcommand{\tensor}{\otimes}
\newcommand{\fol}{\mathcal{F}}
\newcommand{\R}{\mathbb{R}}
\newcommand{\Q}{\mathbb{Q}}
\newcommand{\Z}{\mathbb{Z}}
\newcommand{\N}{\mathbb{N}}
\newcommand{\li}{\mathfrak{l}}
\newcommand{\G}{\Gamma}
\newcommand{\e}{\epsilon}

\begin{document}

\title{The Liouville phenomenon in the deformation problem of coisotropics}
\author{Noah Kieserman}
\thanks{This work was completed with the support of the graduate program of the University of Wisconsin, Madison, a graduate fellowship of the NSF VIGRE program, and the research support of Prof. Yong-Geun Oh.}
\maketitle


\begin{abstract}

The work of Oh and Park ([OP]) on the deformation problem of coisotropic submanifolds opened the possibility of studying a large and interesting class of foliations with some explicit geometric tools. These tools assemble into the structure of an $L_\infty$-algebra on the shifted foliation complex $(\Omega^*[1](\fol), d_\fol)$, which allows a concise description of deformations in terms of a Maurer-Cartan equation. Infinitesimal deformations are given by $d_\fol$-closed forms, and the relation between infinitesimal deformations and full deformations can be studied in terms of obstruction classes lying in the foliation cohomology $H^*_\fol$. Closely related to the foliation cohomology is Haefliger's group $\Omega^*_c(T/H)$, an under-appreciated model for the leaf space of a foliation. We make integral use of this group in showing solvability and unsolvability of the obstruction equations. We also show the $L_\infty$-apparatus to be capable of detecting the Liouville/diophantine distinction of KAM theory, and argue for the greater significance of Haefliger's integration-over-leaves map in passing this fine structure to a geometric model for the leaf space.

\end{abstract}

\tableofcontents

\part{}

\section{Introduction}

Coisotropic submanifolds of symplectic manifolds are a natural generalization of Lagrangian submanifolds. They may equally well be defined on Poisson manifolds, and in that context, as known since [We2], they play precisely the same general role as Lagrangian submanifolds do in the symplectic category: a smooth map $\phi: (M,\pi_M) \to (N,\pi_N)$ between two Poisson manifolds is a Poisson map if and only if its graph is coisotropic in $(M \times N, \pi_M \times -\pi_N)$.

Staying in the symplectic category, a coisotropic is defined as a submanifold $C$ which contains its symplectic orthogonal:
$$(T_xC)^\omega \subseteq T_xC \ \ \forall \ x \in C.$$
This distribution $(TC)^\omega$ on $C$ is of constant rank as $\omega$ is non-degenerate, and is integrable as $\omega$ is closed. Thus a coisotropic has a canonical foliation, the \textit{null foliation}, which we denote $\fol$. In the Lagrangian case it reduces to the trivial foliation by one leaf, and at the other extreme, the symplectic manifold is coisotropic in itself with codimension zero, so the null foliation is the trivial foliation by points. In the general case, we have some mixture of tangential and transverse structure, a foliation whose leaf space is a symplectic space.

The deformation problem of coisotropics in the symplectic category was taken up in 2003 by Oh and Park ([OP]) as a natural boundary condition in the A-model, after Kapustin and Orlov gave evidence that some special coisotropics must be included in the Fukaya category for the homological mirror symmetry conjecture to hold ([KO]). In the same year, Cattaneo and Felder pursued their quantization in Poisson manifolds as natural boundary conditions for the Poisson sigma model ([CF]). The main result of [OP], to be elaborated in Part 1, is:
\begin{theorem}([OP])
There exists the structure of an $L_\infty$-algebra on $\Omega^*[1](\fol)$, the shifted foliation complex of the null foliation. Furthermore, the first-order operator $\li_1$ of this structure is the natural differential $d_\fol$ on this complex, and the Maurer-Cartan equation of this $L_\infty$-algebra is the defining equation for a small deformation of $C$ (the zero section) to be coisotropic.
\end{theorem}

There is a natural obstruction theory for any $L_\infty$-algebra, determining if an infinitesimal deformation may be extended to a full deformation, and thus identifying if the moduli space is singular at a particular point. Here an infinitesimal deformation $\G$ is one for which $d_\fol\G = 0$, and proving that an obstruction class is exact amounts to inverting $d_\fol$. This may be a very badly behaved differential operator, and the associated foliation cohomology is generally infinite-dimensional.

The obstruction classes of an $L_\infty$-algebra are in terms of $\li_2$ and the higher-order operators. In our coisotropic case, $\li_2$ represents the transverse symplectic form. If the foliation has finite holonomy, then the leaf space is a symplectic orbifold, and we may work with the obstruction class directly on the leaf space. It remains to identify a good model for the leaf space in the case of infinite holonomy. A prototypical example is the following: consider a foliation having transverse section a disk, and holonomy given by rotation of the disk by a fixed angle $\alpha$. For $\alpha$ irrational, the naive definition of the leaf space suffers a loss of dimension, and so cannot adequately represent the transverse bracket.

In [Hae], Haefliger offered the group $\Omega^*_c(T/H)$ as a model for the leaf space of a foliation. It is defined in terms of $T$ a transverse section of the foliation, making no connectedness assumption. Here $T$ plays the role of an open cover for the leaf space, with forms on two neighborhoods being identified via pullback along the leaves. This gives a subtler notion than the naive quotient obtained by using forms constant along leaves. Further relevance to our obstruction question is bestowed by Haefliger's integration-over-leaves map $\int_\fol$, which maps forms on the total space to forms in Haefliger's group. Appropriately restricted, this gives a map from foliation cohomology to Haefliger's group.

We will consider an example based on the disk rotation described above, a 4-manifold $Y_\alpha \subset \R^6$, $\alpha \in \R$, viewed as a coisotropic. With the help of Haefliger's map, we find that the behavior of the deformation problem depends not only on if $\alpha$ is irrational, but in fact on what kind of irrational number it is:
\begin{theorem}
The deformation problem of the 4-manifold $Y_\alpha$, as a coisotropic, is obstructed in the two cases
\begin{itemize}
\item $\alpha \in \mathbb{Q}$, and 
\item $\alpha \notin \mathbb{Q},\ \alpha$ a Liouville number, \end{itemize}
and unobstructed in the case \begin{itemize}
\item $\alpha \notin \mathbb{Q},\ \alpha$ satisfying a diophantine condition.\end{itemize}
\end{theorem}

Furthermore, as Haefliger's group already distinguishes the three cases, we view this as deforming some structure over the leaf space, with Haefliger's group as model. We then interpret the result by saying that the deformation problem of the leaf space behaves like the orbifold case ($\alpha \in \mathbb{Q}$), exactly when it is ``close to'' being an orbifold.

\begin{acknowledgements}
I would like to thank Prof. Joel Robbin for many helpful discussions, and especially my advisor, Prof. Yong-Geun Oh, for his generous support during my time in Wisconsin, and for suggesting a question which led to more questions, which led to more questions...
\end{acknowledgements}

\section{$L_\infty$-algebras and their obstruction theory}

We follow the conventions of [DMZ] in our definition of $L_\infty$-algebras, in particular with respect to signs. Let $V$ be a $\Z$-graded vector space, with $v_i$ homogeneous elements of degree $|v_i|$. Denote by $\mathcal{T}(V)$ the tensor algebra of $V$, and by $\Lambda V$ the quotient of $\mathcal{T}(V)$ by the ideal generated by all elements of the form $v_1 \otimes v_2 - (-1)^{|v_1||v_2|}v_2 \otimes v_1$. Then we define the Koszul sign:
\begin{definition}
For $\sigma$ any permutation of $\{1,...k\}$, the \textbf{Koszul sign} $\epsilon(\sigma) \in \{-1,1\}$ is determined by
$$v_1 \land ... \land v_k = \epsilon(\sigma) v_{\sigma(1)} \land ... \land v_{\sigma(k)}.$$
\end{definition}
Note that $\epsilon(\sigma)$ refers to the action of $\sigma$ on $\mathcal{T}(V)$, as it depends on the degrees $|v_i|$: it counts the number of interchanges of odd elements. 
\begin{definition}
We denote the \textbf{antisymmetric Koszul sign} $\chi(\sigma)$ as follows:
$$\chi(\sigma) := sgn(\sigma)\epsilon(\sigma),$$
where $sgn(\sigma)$ refers to the standard sign of $\sigma$ as a permutation of $\{1,...k\}$.
\end{definition}
We call an operator $\mu: \tensor^kV \to V$ \textit{graded commutative} if it is well-defined on $\Lambda V$, that is, if
$$\mu(v_1, ... v_k) = \epsilon(\sigma)\mu(v_{\sigma(1)}, ... v_{\sigma(k)}).$$
We call an operator $\mu$ \textit{graded anti-commutative} if the same identity holds with $\chi$ in place of $\epsilon$.
\begin{remark}
If we consider a Lie algebra as a graded vector space with every element in degree 0, the Lie bracket is graded anti-commutative. If we consider every element to be in degree 1, the bracket is graded commutative.
\end{remark}
Now for the definition of $L_\infty$-algebra, we wish to form the composition of two graded anti-commutative operators.
\begin{definition} Given two graded anti-commutative operators $\mu : \otimes^kV \to V$ and $\nu : \otimes^jV \to V$, we define the graded anti-commutative operator $\mu \circ \nu : \otimes^{j+k-1}V \to V$, by
$$\mu \circ \nu (v_1, ...  v_{j+k-1}) := \sum_{\substack{\sigma \in S_{j,k-1}}} (-1)^{deg(\nu)} \chi(\sigma) \mu(\nu(v_{\sigma(1)},... v_{\sigma(j)}), v_{\sigma(j+1)}, ... v_{\sigma(j+k-1)}).$$
Here $S_{j,k-1}$ denotes the subgroup of $S_{j+k-1}$ consisting of all unshuffles of $j+k-1$ elements into two ordered sets of $j$ elements and $k-1$ elements, respectively.
\end{definition}

Finally we are ready to define
\begin{definition}
An \textbf{$L_\infty$-algebra} is $V$, a $\Z$-graded vector space, together with a sequence of graded anti-commutative operators
$$\li_k : \otimes^kV \to V,\ \ k \in \N,$$
of degree 2-k, such that
$$\sum_{\substack{n}} \sum_{\substack{i+j=n}} \li_i \circ \li_j = 0.$$
\end{definition}

\begin{remark}
Standard terminology designates as \textit{weak} an $L_\infty$-algebra in which $\li_0 \neq 0$. The condition $\li_0=0$ is also referred to as \textit{flatness}. As noted in [OP], the $L_\infty$-algebra governing the deformation of coisotropics is flat, and we assume flatness for the remainder of the present work. Some light is shed on flatness in the Poisson context by the approach of Voronov via higher derived brackets ([Vo]). Cattaneo and Felder, in [CF], show that the $L_\infty$-structure produced for any submanifold of a Poisson manifold is flat if and only if the submanifold is coisotropic. For a further study of higher homotopy structures and their application to the deformation of coisotropics, see [Sch].
\end{remark}

For a flat $L_\infty$-algebra, the first coherence relations take the form
$$\begin{array}{lcl}
\li_1 \circ \li_1 &=& 0\\
\li_1 \circ \li_2 + \li_2 \circ \li_1 &=& 0\\
\li_1 \circ \li_3 + \li_2 \circ \li_2 + \li_3 \circ \li_1 &=& 0\\
&\vdots&\\
\end{array}$$
Thus the first relation allows us to define a cohomology theory, and the second encodes the Leibniz rule for $\li_1$ and $\li_2$, almost yielding a differential graded Lie algebra (dgLa). The third relation measures the failure of the Jacobi identity $(\li_2 \circ \li_2 \neq 0)$, a phenomenon we will return to in our coisotropic case below. 
\begin{definition}
For $(V,\{\li_k\})$ an $L_\infty$-algebra and $v \in V$, the \textit{Maurer-Cartan equation} (MC) is
$$\sum \frac{1}{k!}\li_k(v, ... v) = 0,$$
which reduces to the traditional case $dv + \frac{1}{2}\{v,v\} = 0$ if $\li_k=0\ \forall\ k \geq 3$.
\end{definition}
An $L_\infty$-algebra is said to \textit{govern} a deformation problem if deformations preserving the structure under consideration are exactly those satisfying the MC equation of the given $L_\infty$-algebra.

To study deformations using this apparatus, we attempt to find power series solutions to the MC equation near zero. That is, 
$$\Gamma_t = \sum_{\substack{i=1}}^{\substack{\infty}} \Gamma_i t^i,$$
where $t$ is a formal parameter, and each $\Gamma_i$ is an element of the underlying vector space $V$. For our purposes it suffices to consider degree 1 elements (shifted degree 0), incidentally simplifying signs considerably. We require $\Gamma_t$ to solve the MC equation for each $t$, and after gathering terms by degree, we find:
$$\begin{array}{rcl}
t &:& d \Gamma_1 = 0;\\
\\
t^2 &:& -d \Gamma_2 = \frac{1}{2}\{\Gamma_1,\Gamma_1\};\\
\\
t^3 &:& -d \Gamma_3 = \frac{1}{3!}\li_3(\Gamma_1, \Gamma_1, \Gamma_1) + \frac{1}{2}(\{\Gamma_1,\Gamma_2\} + \{\Gamma_2,\Gamma_1\});\\
\\
t^4 &:& -d\G_4 = \frac{1}{2}(\{\G_3,\G_1\} + \{\G_2, \G_2\} + \{\G_1,\G_3\})\\ &\ &\ \ \ \ \ \ \ \ \ \ \  + \frac{1}{3!}(\li_3(\G_2, \G_1, \G_1) + \li_3(\G_1, \G_2, \G_1) + \li_3(\G_1, \G_1, \G_2))\\ &\ &\ \ \ \ \ \ \ \ \ \ \ + \frac{1}{4!}\li_4(\G_1, \G_1, \G_1, \G_1)\\
& &\ \ \  \vdots
\end{array}$$

We call $\Gamma_1$ satisfying $d \Gamma_1 = 0$ an \textit{infinitesimal deformation}. From the $L_\infty$-relations, it follows that $d\{\G_1,\G_1\} = \{d\G_1,\G_1\} + \{\G_1, d\G_1\} = 0$. Then the second equation above says that the $\li_1$-cohomology class $[\{\Gamma_1,\Gamma_1\}]$ vanishes, if and only if there exists $\Gamma_2$ extending $\G$ to second order. In general, for any $i$, given $\G_1,...\G_{i-1}$ satisfying these equations up to order $i-1$, the right hand side of equation $i$ defines a cohomology class, the $i^{th}$ obstruction class, which vanishes iff $\G$ extends to order $i$. 

The first obstruction class $[\{\G_1,\G_1\}]$, sometimes referred to as the value of the \textit{Kuranishi map} $[\G_1] \mapsto [\{\G_1,\G_1\}]$, suffices to prove our obstructedness results. We return to the higher obstructions in the final unobstructedness section.

\subsection{Remark on the minimal model}

As $\li_1$ (hereafter referred to as $d$) and $\li_2$ satisfy the Leibniz rule, $\li_2$ induces an operation on $\li_1$-cohomology. In addition, if $\G_i$ are closed for $d=\li_1$,
$$\begin{array}{rcl}
\li_2 \circ \li_2 (\G_1,\G_2,\G_3) &=& - (d \circ \li_3 + \li_3 \circ d) (\G_1,\G_2,\G_3)\\
Jacobi (\li_2;\G_1,\G_2,\G_3) &=& - d\li_3(\G_1,\G_2,\G_3) \pm \sum_{\substack{\sigma}} \li_3(d\G_{\sigma(1)},\G_{\sigma(2)},\G_{\sigma(3)})\\
Jacobi (\li_2;\G_1,\G_2,\G_3) &=& d\Sigma + 0 = d\Sigma,\\
\end{array}$$
for some $\Sigma \in V$. Thus $\li_2$ in fact defines an honest Lie bracket on $\li_1$-cohomology, and $(H^*_{\li_1}, d=0, [\li_2])$ gives a dgLa.

On the other hand, rewriting the above, we see $\li_3$ does not necessarily descend to cohomology, much less as zero:
$$\begin{array}{rcl}
d \circ \li_3(\G_1,\G_2,\G_3) &=& - \li_3 \circ d (\G_1,\G_2,\G_3) - Jacobi (\li_2;\G_1,\G_2,\G_3)\\
&=& 0 - d\Sigma\\
\end{array}$$
($\li_3$ of three closed forms is not necessarily closed). Therefore the question remains of what relationship there is between the original $L_\infty$-algebra, and the induced $L_\infty$ (dgLa) structure on its $\li_1$-cohomology. While we do not answer this question in the present work, we hope to shed some light on its meaning for coisotropics.

\section{Oh and Park's $L_\infty$-algebra for deformation of coisotropics}

In this section we give a tensor-level introduction to Oh and Park's $L_\infty$-structure, as identified in [OP]. Recall:
\begin{definition}
Let $(M,\omega)$ be a symplectic manifold. A submanifold $C \subseteq M$ is called \textbf{coisotropic} iff
$$(T_xC)^\omega \subset T_xC\ \ \forall x \in C,$$
where $(\cdot)^\omega$ denotes the symplectic orthogonal.
\end{definition}
The assumption of closedness of $\omega$ passes to $C$. Thus $\omega|_C$ generates a differential ideal in $\Omega^*(C)$, which is equivalent to integrability of $ker(\omega|_C)$ by the Frobenius Theorem. We call this, the foliation $\fol$ of $C$ such that $T\fol = ker(\omega|_C)$, the \textit{null foliation} of $(C, \omega|_C)$.

There is also what remains of the symplectic form, "on the quotient $C/\fol$." To be more precise, there is a non-degenerate form on $N\fol := TC/T\fol$, which in the case of $\fol$ given by a fiber bundle, may be interpreted as a symplectic form on the base. The deRham differential does not canonically restrict to $N^*\fol$, however, and in general we only have a non-degenerate form, not necessarily closed. (On the other hand, it would be a misnomer to call it an \textit{almost symplectic structure}, for the failure of closedness is not a defect of the form itself, but rather of the underlying space in its failure to be a manifold.)

Thus there are two aspects to the structure of a coisotropic, the null foliation and the transverse structure, and both present obstacles when considering deformations. The deformation of foliations is quite general and involved, and while the deformation of (compact) symplectic manifolds is unobstructed, we rapidly encounter leaf spaces which are not manifolds. It would represent a beginning toward both problems to identify a fine model for the leaf space. 
\vspace{.2in}

One model for the leaf space of a foliation is the \textit{basic} forms.
\begin{definition}
Let $(M,\fol)$ be a smooth, foliated manifold. A differential form $\alpha \in \Omega^*(M)$ is called \textbf{basic} if, for all $X \in T\fol$,
$$\iota_X \rfloor \alpha = \iota_X \rfloor d\alpha = 0$$
In particular, a function $f \in C^\infty(M)$ is basic if for all $X \in T\fol$, $\mathcal{L}_Xf = 0$.
\end{definition}
Proposition 5.2 of [OP] states holonomy-invariance of $\omega|_C$, which, combined with closedness, implies it is basic. In our example, the basic cohomology will be rather easy to describe, but carries very little information.
\vspace{.15in}

Another tool, which keeps track of essentially the opposite information, is the \textit{foliation cohomology} $H^*_\fol$, defined from the cochain complex $(\Omega^*_\fol, d_\fol)$, which we describe in the next section. The principal problem with this complex is that unlike basic forms, it does not canonically include in $\Omega^*(C)$. There is a spectral sequence decomposing $\Omega^*(C)$ ([To]), and given an inclusion $\Omega^*_\fol(C) \hookrightarrow \Omega^*(C)$, the foliation cohomology $H^*_\fol(C)$ appears in the first column at the $E_1$ level, the complex of basic forms $\Omega^*_{bas}(C)$ in the first row. Hence $H^0_\fol(C) \cong \Omega^0_{bas}(C) \equiv C^\infty_{bas}(C)$, and all of $H^*_\fol$ is a module over basic functions. Even basic cohomology may already be infinite-dimensional ([Gh]), so at first glance foliation cohomology is a rather poorly behaved object. What will concern us is exactly the part of $H^*_\fol$ that goes beyond $\Omega^0_{bas}$.
\vspace{.2in}

Finally, for completeness we mention two facts. One is that were we to restrict attention to transversely holomorphic foliations, there is a versality theorem for transversely symplectic foliations that would apply, due to Girbau and Guasp ([GG]). The other is that there is a \v{C}ech-deRham theory for foliations, defined by Crainic and Moerdijk ([CM]), which subsumes the above cohomologies and has better properties. Among other things, it satisfies Poincar\'{e} duality. There is a natural map to foliation cohomology with compact supports, and in general some information is lost. The foliation \v{C}ech-deRham theory, however, is difficult to compute with.

\subsection{The cochain complex $(\Omega^*(\fol), d_\fol)$}

In the Lagrangian embedding theorem ([We1]), a tubular neighborhood of a Lagrangian submanifold $L$ of a symplectic manifold is canonically identified with the cotangent bundle of $L$ with the standard symplectic form. Furthermore, ($C^1$-small) sections are Lagrangian if and only if they are closed as 1-forms. Thus the deformation problem of Lagrangian submanifolds is unobstructed, in fact linear. After identifying trivial deformations with exact forms, we find it has a finite-dimensional moduli space, modelled on $H^1(L;\R)$. 

Correspondingly, by a coisotropic embedding theorem of Gotay [Go], the normal bundle of $C$ in $M$ is modelled (non-canonically) by the cotangent bundle $T^*\fol$ of the null foliation. From this one constructs the foliation deRham complex $\Omega^*(\fol) := \G (\Lambda^*(T^*\fol))$, which hosts a canonical differential $d_\fol$ for any foliation $\fol$. This differential is defined by the Cartan formula, which makes sense due to involutivity of $T\fol$:
\begin{definition}
Let $\alpha \in \Omega^k(\fol) := \G(\Lambda^k T^*\fol)$, and $X_i \in T\fol$. Then
$$\begin{array}{rcl}
d_\fol\alpha(X_0, ... X_k) &:=& \sum (-1)^i X_i(\alpha(X_0,...\hat{X_i}, ... X_k)) \\
& &\ \ \ \ + \sum (-1)^{i+j} \alpha([X_i,X_j],X_0, ... \hat{X_i}, ... \hat{X_j}, ... X_k).
\end{array}$$
\end{definition}
Thus $H^*_\fol := H^*(\Omega^*\fol, d_\fol)$ may be viewed alternatively as Lie algebra cohomology of the foliation, Lie algebroid cohomology of the tangent bundle of the foliation, and the restriction of deRham cohomology to the foliation. (For further discussion of foliated cohomology and references, see [CM], [dS], [Va].) The first way in which the coisotropic case differs from the Lagrangian case is that this cohomology is generally infinite dimensional. 

There is a symplectic form on $T^*\fol$ as well ([Go]), again non-canonical, which allows us to ask when the graph of a small section is coisotropic. Coisotropic deformations also differ from Lagrangians in that they are not described entirely by $d_\fol$. Rather, the defining equation is non-linear to all orders, given by the Maurer-Cartan equation of an $L_\infty$-algebra. The higher-order operators of the $L_\infty$-algebra include the transverse symplectic structure, and additionally, the defect in the symplectic reduction as a smooth manifold. 

 The first-order operator $\li_1$ of this $L_\infty$-structure is given by $\li_1(\G) := (-1)^{|\G|}d_\fol\G$. We now proceed to define $\li_2$.

\subsection{Transverse Poisson bracket}

To define higher-order operators on $\Omega^*(\fol)$, we note that $\omega$ induces a well-defined non-degenerate form on $N\fol := TC/T\fol$. It is not possible to speak of closedness of this form, however, as the bundle $N^*\fol$ itself, while a subbundle of $T^*C$, is not closed under $d|_C$. In fact we consider the Poisson bivector dual to the symplectic form on $T^*\fol$. The non-degenerate form on $N\fol$ also has a dual $P \in \Lambda^2 N\fol$, not necessarily satisfying the Jacobi identity. To write $P$ in coordinates, we must choose a representative for the normal bundle. That is, we choose $G$ a subbundle of  $TC$ such that $\forall x \in C$,
$$T_xC \cong G_x \oplus T_x\fol,$$
otherwise known as an \textbf{Ehresmann connection} for the foliation. For any choice of local foliation coordinates $(y^i, q^\alpha)$, i.e. local coordinates on $C$ such that $\fol$ is defined by $y^i = const$, we write $G = span\{e_i\}$, where
$$e_i := \frac{\partial}{\partial y^i} + R^\alpha_i(y^i,q^\alpha) \frac{\partial}{\partial q^\alpha}.$$
(It is equivalent to make an arbitrary choice of $R \in \G(T^*C \otimes T\fol)$, where $R(\frac{\partial}{\partial y^i}) = -R^\alpha_i \frac{\partial}{\partial q^\alpha}$.) Then we may write
$$P = P^{ij}e_i \land e_j,$$
where $(P^{ij}) = (\omega_{ij})^{-1},$ and $\omega_{ij} = \omega(e_i, e_j) = \omega(\frac{\partial}{\partial y^i}, \frac{\partial}{\partial y^j})$.

Such a choice also allows us to make sense of taking transverse derivatives. A lengthy but straightforward derivation in Section 7 of [OP] shows that
\begin{enumerate}
\item The choice of $G$ induces a connection $\nabla$ on the bundle $T^*\fol$, that is,
$$(\nabla\Gamma)_x \in \Gamma(T^*_xC \tensor T^*_x\fol).$$
\item In coordinates,
$$\nabla \G  = \frac{\partial\G_\alpha}{\partial q^\beta}\ dq^\beta \otimes (\frac{\partial}{\partial q^\alpha})^*\ \oplus\  (\frac{\partial\G_\alpha}{\partial y^i} + \G_\beta \frac{\partial R^\beta_i}{\partial q^\alpha})\ dy^i \otimes (\frac{\partial}{\partial q^\alpha})^*.$$
\end{enumerate}
Given $G \cong N\fol$, we may write $T^*\fol$ as $G^\circ$, the annihilator of $G$, with basis written in coordinates as $(\frac{\partial}{\partial q^\alpha})^* = dq^\alpha - R^\alpha_idy^i$. Then 
$$\nabla \G  = \nabla_\beta \G (\frac{\partial}{\partial q^\beta})^* + \nabla_i \G (e_i)^*,$$
where we denote
$$\begin{array}{rcl}
\nabla_\beta \G  &=& \nabla_\beta\Gamma_\alpha(\frac{\partial}{\partial q^\alpha})^* = \frac{\partial\G_\alpha}{\partial q^\beta} (\frac{\partial}{\partial q^\alpha})^*\\
\nabla_i \G &=& \nabla_i \G_\alpha (\frac{\partial}{\partial q^\alpha})^* = (\frac{\partial\G_\alpha}{\partial y^i} + R^\beta_i\frac{\partial \G_\alpha}{\partial q^\beta} + \G_\beta \frac{\partial R^\beta_i}{\partial q^\alpha}) (\frac{\partial}{\partial q^\alpha})^*
\end{array}$$
(Note that in [OP], this notation is used to signify an earlier stage.) Here $d_\fol$ may be seen in yet another way, as the skew-symmetrization of $\nabla_\beta$: since $\nabla$ is flat in the $T\fol$ direction, skew-symmetrizing $\nabla_\beta$ yields a differential.

Denote by $\Omega^*[1](\fol)$ the degree shift defined by  $\Omega[1]^k(\fol) := \Omega^{k+1}(\fol)$. We may now state:
\begin{theorem} ([OP], Theorem 9.4)
Given a coisotropic $C$ and a choice of transverse distribution $G$ complementary to the null foliation, there exists an $L_\infty$-structure on $\Omega[1]^*(\fol)$, such that for $\Sigma \in \Omega[1]^k(\fol)$ and $\G_i \in \Omega[1]^0(\fol)$, 
$$\begin{array}{rcl}
\li_1(\sigma) &=& (-1)^k d_\fol\Sigma\\
\li_2(\G_1,\G_2) &=& P(\nabla \G_1,\nabla \G_2).
\end{array}$$
\qed
\end{theorem}
\begin{theorem}([OP], Theorem 10.1)
The $L_\infty$-structures induced by different choices of transverse distribution $G$ are canonically isomorphic, so it is an invariant of the coisotropic. \qed
\end{theorem}

\begin{remark}($P_\infty$-structure)
The vector space $\Omega^*(\fol)$ of the $L_\infty$-structure for the deformation problem of coisotropics is in fact an algebra. We note that the $L_\infty$-operators $\li_1$ and $\li_2$ defined thus far have an additional property, namely that they are multiderivations for the multiplication of sections. This global nature is recognized by referring to an "$L_\infty$-algebroid" in [OP]. In [CF] it is stated explicitly, defining a "$P_\infty$-structure" to be an $L_\infty$-algebra for which the underlying vector space is an algebra, and all of the operators $\li_k$ are multiderivations for the product.
\end{remark}
\vspace{.2in}

\subsection{Obstruction}

Consider an infinitesimal deformation, that is, $\G \in \G(T^*\fol)$ such that $d_\fol \G = 0$. We wish to prove exactness of obstructions, i.e. existence of $\Sigma \in \G(T^*\fol)$ such that $d_\fol\Sigma = \{\G,\G\}$. As it will be relevant in our main example, we consider the case of a dimension-2 foliation, where we have, for $h_1, h_2, f,$ and $g$ smooth functions on $C$,
$$\begin{array}{rcl}
\Sigma &=& h_1(\frac{\partial}{\partial q^1})^* + h_2(\frac{\partial}{\partial q^2})^*,\\
\Gamma &=& f(\frac{\partial}{\partial q^1})^* + g(\frac{\partial}{\partial q^2})^*.\\
\end{array}$$
Then $d_\fol$ is computed:
$$\begin{array}{rcl}
d_\fol\Sigma &=& d_\fol (h_1(\frac{\partial}{\partial q^1})^* + h_2(\frac{\partial}{\partial q^2})^*)\\
&=& d_\fol(h_1) \land (\frac{\partial}{\partial q^1})^* + d_\fol(h_2) (\frac{\partial}{\partial q^2})^*\\
&=& \frac{\partial h_1}{\partial q^2}(\frac{\partial}{\partial q^2})^* \land (\frac{\partial}{\partial q^1})^* + \frac{\partial h_2}{\partial q^1} (\frac{\partial}{\partial q^1})^* \land (\frac{\partial}{\partial q^2})^*\\
&=& (\frac{\partial h_2}{\partial q^1} - \frac{\partial h_1}{\partial q^2})\  (\frac{\partial}{\partial q^1})^* \land (\frac{\partial}{\partial q^2})^*.\\
\end{array}$$
The first obstruction is given at the chain level by
$$\begin{array}{rcl}
\li_2(\Gamma,\Gamma) &:=& P(\nabla\Gamma,\nabla\Gamma)\\
&=&P^{ij}\ \nabla_i(\Gamma) \land \nabla_j(\Gamma)\\
&=&P^{ij}\ (\nabla_i(f)\nabla_j(g) - \nabla_j(f)\nabla_i(g))\ (\frac{\partial}{\partial q^1})^* \land (\frac{\partial}{\partial q^2})^*,
\end{array} $$
where there is summation over $i,j$ all transverse coordinates. Here we introduce the following notation
\begin{definition}
For $C$ a coisotropic with transverse bivector $P$, $\nabla$ a choice of connection on $T^*\fol$, and $f,g \in C^\infty(C)$,
$$\{f,g\}_P := P^{ij}\ (\nabla_i(f)\nabla_j(g) - \nabla_j(f)\nabla_i(g))$$
\end{definition}
This clearly depends on the choice of $\nabla$ (although Theorem 10.1 of [OP] tells us that the cohomology class of $\li_2(\Gamma,\Gamma)$ does not). In our obstructed cases below, however, it will suffice to work with functions which can be written independently of the tangential coordinates. What's more, it will suffice to work on an open set where $\nabla$ is flat, meaning it specifies coordinates for a transverse submanifold $T$. Then $P|_T$ makes an honest bracket on $T$, and $\{f,g\}_P$ will have independent meaning of the choice of $T$.

With the above understanding of notation, our main obstructedness equation becomes
\begin{equation}
\begin{array}{rcl}
d_\fol \Sigma &=& \{\G,\G\}\\
\frac{\partial h_2}{\partial q^1} - \frac{\partial h_1}{\partial q^2} &=& \{f,g\}_P.\\
\end{array}
\end{equation}
\vspace{.1in}
\begin{example}
If we consider $\R^6$ with the standard symplectic form $\sum_{\substack{i=1}}^{\substack{3}}dq_i \land dp^i$, and
$$C = \R^4 = \{(q_1,p^1,q_2,0,q_3,0)\}$$
embedded as a coisotropic submanifold, the null foliation has leaves $\R^2$ parametrized by $(q_2,q_3)$, and leaf space $\R^2$ parametrized by $(q_1,p^1)$ with the standard symplectic form on $\R^2$.
$$\begin{array}{rcl}
d_\fol \Sigma &=& \{\G,\G\}\\
\frac{\partial h_2}{\partial q^2} - \frac{\partial h_1}{\partial q^3} &=& \frac{\partial f}{\partial q^1}\frac{\partial g}{\partial p^1} - \frac{\partial g}{\partial q^1}\frac{\partial f}{\partial p^1}.\\
\end{array}$$
We may always solve for $\Sigma$, integrating along $q_2$ and $q_3$. In fact, the higher $L_\infty$-operators vanish as well, and this deformation problem is unobstructed.
\end{example}
\vspace{.1in}
\begin{example}
In an example considered from an analytic point of view by Zambon [Za], the torus $T^4 = \{(\theta_1,\theta_2,\theta_3,\theta_4)\}$ is considered as a pre-symplectic manifold (an equivalent problem), with presymplectic form $d\theta_1 \land d\theta_2$. The leaves of the null foliation are tori $T^2$ parametrized by $(\theta_3,\theta_4)$, and the leaf space is $T^2$ parametrized by $(\theta_1,\theta_2)$. There is a non-trivial bracket on the leaf space, and furthermore, there is a precondition for solution of the obstruction equation $d_\fol\Sigma = \{\G,\G\}$, namely of zero average over the leaves. This is not in general satisfied, and thus due to the existence of leaf topology, this example is obstructed.
\end{example}
\vspace{.1in}
Thus obstructedness arises from a non-trivial transverse bracket, together with leaf topology as reflected in the operator $d_\fol$, specifically in its invertibility. The above foliations are trivial - below we will define our main example with the goal of examining the effects of a non-trivial foliation on obstructedness.
\vspace{.1in}

\subsection{Transverse curvature $F_G$}

The first obstruction will suffice to prove our obstruction results, and we postpone the full definition of the higher $L_\infty$-operators to the section on unobstructedness below. Here we simply clarify the statements made above that the conormal bundle $N^*\fol$ is not closed under $d|_C$, and that the transverse distribution $G$ is not always integrable.

The distribution $G$ has a well-defined curvature, which vanishes iff $G$ is integrable. The definition is a standard one for a connection given as a vector-valued 1-form ([Sa]). The fundamental object there is the connection
$$R \in \G(T^*C \otimes T\fol),$$
where $G = ker(R)$, viewing $R$ as a projection $TC \to T\fol$. This projection is referred to as $\Pi_G$ in [OP]. Then we may define the curvature by:
\begin{definition}
For $\eta_i = [X_i] \in \G(N\fol),\ X_i \in G$ (which is unique), we define $F_G \in \G(\Lambda^2N^*\fol \otimes T\fol)$ by
$$F_G(\eta_1, \eta_2) := R([X_1,X_2]).$$
\end{definition}
This is referred to as $F_\Pi$ in [OP], and is a direct restriction of [Sa], p. 89 from $T^*C$ to $N^*\fol$. 
\begin{proposition}([OP], Proposition 4.2) 
This is tensorial. \qed
\end{proposition}
By a computation, we find that for
$$F_G = F^\alpha_{ij} \frac{\partial}{\partial q^\alpha} \tensor (e_i)^* \land (e_j)^*,$$
$$F^\alpha_{ij} = \left(\frac{\partial R^\alpha_j}{\partial y^i} + R^\gamma_i\frac{\partial R^\alpha_j}{\partial q^\gamma} \right)_{skew},$$
where "skew" refers to the indices $i,j$. This expression will be used in the final unobstructedness section. The expression for the higher-order operators $\li_k,\ k \geq 3$, incorporates the curvature, and the higher operators vanish, yielding a dgLa, if and only if the curvature is zero. On the other hand, the curvature vanishes if and only if $G$ integrates to give a transverse foliation, complementary to the null foliation. We may attempt to improve our choice of $G$ and stay within the same $L_\infty$-isomorphism class, yet there exist foliations for which no transverse foliation exists, that is, for which no choice of $G$ can be integrable. We end this expository part with the purely algebraic question:

\begin{question}
When does a given $L_\infty$-algebra have a dgLa representative in its $L_\infty$-isomorphism class?
\end{question}
\vspace{.2in}




\part{}

\section{Definition of the main example}

Consider $\R^{2n}$ with the standard symplectic form $\omega_0$. We use polar coordinates $(r_i, \theta_i), i = 1, \ldots n$, and the notation $H_i := \frac{1}{2}r_i^2$. We define, for $\alpha \in \R, \alpha > 1,$
$$H_\alpha := H_1 + \alpha H_2 : \R^4 \to \R.$$
We will abuse notation and also refer to $H_\alpha$ on $\R^6$. Then define
$$M_\alpha := H_\alpha^{-1}(\frac{1}{4}) \subset \mathbb{R}^4,$$
which in $\R^6$ yields $M_\alpha \times \R^2$. The coisotropic whose deformation problem we consider is defined as follows:
$$\begin{array}{rcl}
Y_\alpha &:=& H_\alpha^{-1}(\frac{1}{4}) \cap  H_3^{-1}(\frac{1}{2}) \subset \R^6\\
&\cong& M_\alpha \times S^1.\\
\end{array}$$ 
Thus the definitions of $M_\alpha$ and $Y_\alpha$ determine either one of $r_1, r_2$ in terms of the other, and lead to complementary bounds:
$$\begin{array}{rcl}
r_1 &\in& [0, \frac{1}{\sqrt{2}}]\\
r_2 &\in& [0, \frac{1}{\sqrt{2\alpha}}].
\end{array}$$

As a submanifold of the symplectic manifold $(\R^6,\omega_0)$, $Y_\alpha$ is a coisotropic submanifold. As the vector fields $X_{H_\alpha}$ and $X_{H_3}$ come from commuting Hamiltonians, they are both tangent to their common level set, $Y_\alpha$. Furthermore, being that their flow preserves the symplectic form on $\R^6$, they are in the kernel of $\omega_0|_{Y_\alpha}$. Writing them explicitly, 
$$\begin{array}{rcl}
X_{H_\alpha} &=& \frac{\partial}{\partial \theta_1} + \alpha \frac{\partial}{\partial \theta_2}\\
X_{H_3} &=& \frac{\partial}{\partial \theta_3},\\
\end{array}$$
we see they are in addition everywhere nonzero on $Y_\alpha$ and linearly independent, so they span all of the maximally 2-dimensional kernel. We denote this, the null foliation, by $\fol$, so $T\fol = span\{X_{H_\alpha}, X_{H_3}\}$. Likewise, and more simply, $M_\alpha$ is coisotropic in $\R^4$ with 1-dimensional null foliation $\tilde{\fol}$ spanned by $X_{H_\alpha}$. These two, $(M_\alpha,\tilde{\fol})$ and $(Y_\alpha,\fol)$, yield homeomorphic leaf spaces, and a complete transversal $T$ to $\tilde{\fol}$ yields a complete transversal $T \times \{*\}$ to $\fol$. (A complete transversal is a submanifold transverse to the foliation, not connected in general, which intersects each leaf at least once.) For the moment we restrict our attention to $M_\alpha$.

 $M_\alpha$ is simply an ellipsoid in $\R^4$, i.e. diffeomorphic to $S^3$. We note that the radial coordinates in $\R^4$ are holonomy-invariant, so the foliation is confined to tori of fixed radii, except when one of the radial coordinates is zero, where the flow of $X_{H_\alpha}$ is along one of the two complementary center circles $r_2 = 0$ and $r_1 = 0$. We have obtained the picture of $M_\alpha$ as $S^3$ decomposed into complementary families of nested tori $S^1(r_1) \times S^1(r_2)$:
\begin{center}
\includegraphics[height=3in]{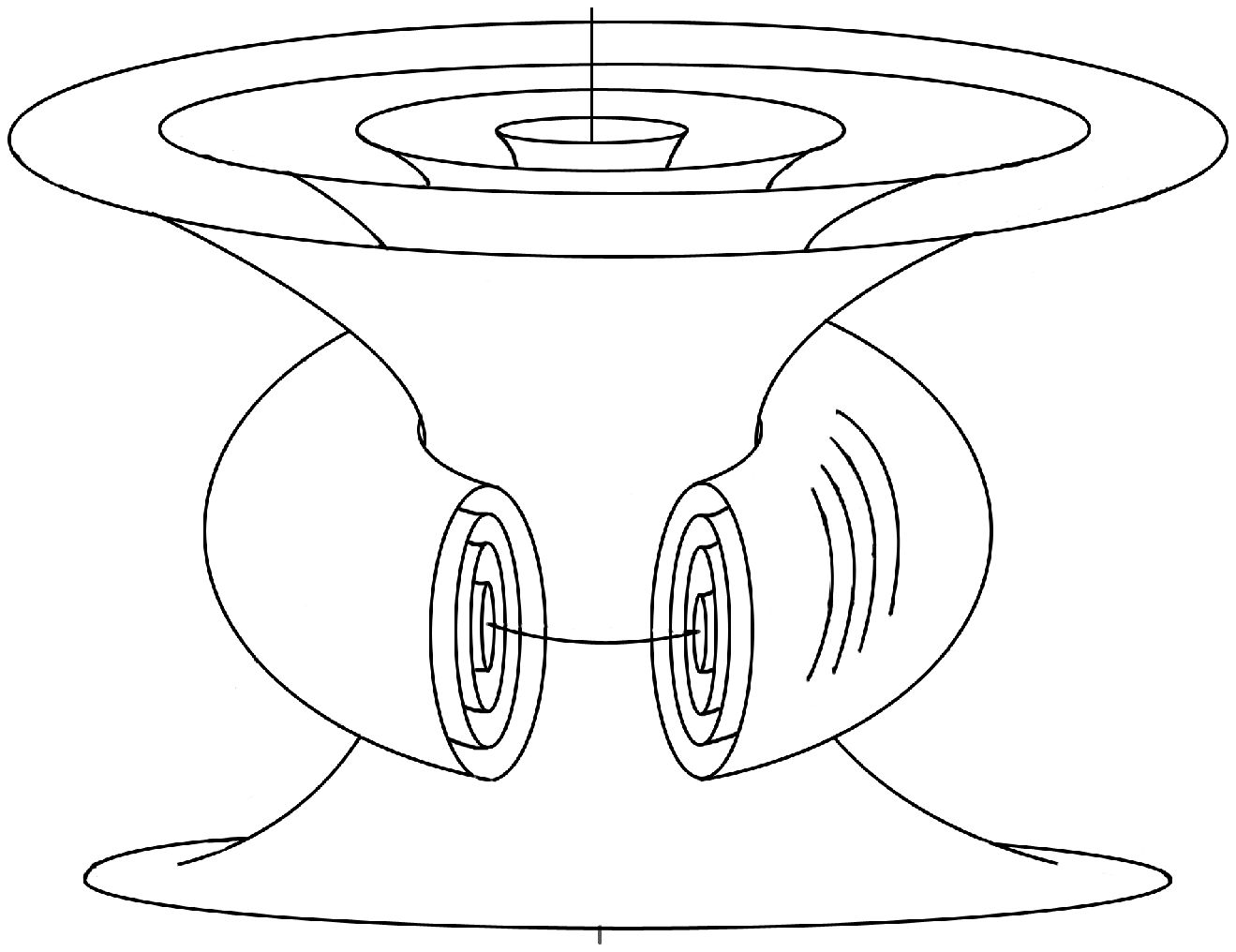}
\end{center}
$$M^3\ \cong\ S^3\ \cong\ S^1 \times D^2\ \amalg_{T^2}\ D^2 \times S^1.$$
Each torus is foliated by lines of slope $\alpha$, and there are in addition two exceptional leaves, the center circles: $S^1(r_1) \times \{0\}$ and $\{0\} \times S^1(r_2)$. In particular we see that a complete transversal is given by two disks $D^2 \times \{*\}\ \coprod\ \{*\} \times D^2$. Returning to $Y_\alpha$, $X_{H_3}$ carries no holonomy, and the foliation $\fol$ is exactly $\tilde{\fol} \times S^1$.

In the case $\alpha \in \Q$, all leaves are closed, diffeomorphic to $T^2$. Only the two exceptional leaves have non-trivial holonomy: if $\alpha = p/q$, one leaf has holonomy $\Z/p\Z$ and the other $\Z/q\Z$. The leaf space is then an orbifold, topologically $S^2$, covered by the above transverse disks, and having two orbifold points at the ``poles." Its deformation will be seen to be obstructed for the same reason as before: leaf topology giving rise to a zero average condition over the leaves. The main question is what happens if $\alpha \notin \Q$. In fact it will depend on what sort of irrational number $\alpha$ is, namely if it is a Liouville number or not. Our result is the following, the proof of which constitutes the remainder of the paper.
\begin{theorem}
The deformation problem of the 4-manifold $Y_\alpha$, as a coisotropic, is obstructed in the two cases
\begin{itemize}
\item $\alpha \in \mathbb{Q}$, and 
\item $\alpha \notin \mathbb{Q},\ \alpha$ a Liouville number, \end{itemize}
and unobstructed in the case \begin{itemize}
\item $\alpha \notin \mathbb{Q},\ \alpha$ satisfying a diophantine condition.\end{itemize}
\end{theorem}
\vspace{.2in}

\section{Haefliger's integration-over-leaves map}

To show an infinitesimal deformation $[\Gamma] \in H^1_\fol(Y_\alpha)$ to be obstructed is to show non-exactness of its first obstruction class $\{\G,\G\}$ in $H^2_\fol(Y_\alpha)$. Rather, we will consider its image under an isomorphism, Haefliger's map $\int_\fol$, and show that to be non-zero.

We recall the definition of Haefliger's group [Hae], for differential forms of any degree. By a \textit{complete transversal} we mean an immersed submanifold everywhere transverse to the foliation, that intersects each leaf at least once. Of note is the lack of connectedness assumption.
\begin{definition}
Let $T$ be a complete transversal for a foliation defined by the holonomy pseudogroup $H$. Then denoting $h \in H$ (so $h$ is a locally-defined diffeomorphism of $T$), and $\tau \in \Omega^k_c(T)$ having support in the range of $h$, we define the group $\Omega^k_c(T/H)$ as the vector space quotient
$$\Omega^k_c(T/H) := \Omega^k_c(T)/span_\R\{\tau - h^*\tau\}.$$
\end{definition}
With a suitable notion of equivalence of pseudogroups, this definition is independent of the choices of $T$ and $H$ ([Hae] \S1.2).

Haefliger then defines the map $\int_\fol$, ``integration-over-leaves,'' defined at the chain level:
$$\int_\fol : \Omega_c^{p+k}(M) \to \Omega_c^k(T/H),$$
where $p = dim(\fol)$. The general case proceeds as follows. One first takes a \textit{regular cover} $\mathcal{U}$ of $(M,\fol)$: 
\begin{definition}
For $\mathcal{U} = \{U_\alpha\}$ an open cover of $(M,\fol)$, we denote by $\pi_\alpha:U_\alpha \to T_\alpha := U_\alpha/ \fol$ the quotient map onto the local leaf space. We say the cover $\mathcal{U}$ is $\mathbf{regular}$ if
\begin{enumerate}
\item $\pi_\alpha$ is a submersion onto $T_\alpha$ a smooth manifold, and
\item each plaque in $U_\alpha$ meets at most one plaque in $U_\beta$.
\end{enumerate}
(A plaque is a connected component of $\pi_\alpha^{-1}(x)$.)
\end{definition}
\begin{example}
Consider the foliation of $T^2 = \{(\theta_1,\theta_2)\}$ given by lines $\theta_1 = const$. Then a regular cover is given by three opens, $\theta_2 \in (-\e, \frac{2\pi}{3}+\e), (\frac{2\pi}{3}-\e, \frac{4\pi}{3}+\e), (\frac{4\pi}{3}-\e, \frac{2\pi}+\e)$, with complete transversal three circles, as pictured below. Here the regular cover is clearly unnecessary, as there is no holonomy to interfere with well-definedness of a global integration-over-fibers.
\begin{center}
\includegraphics[width=5in]{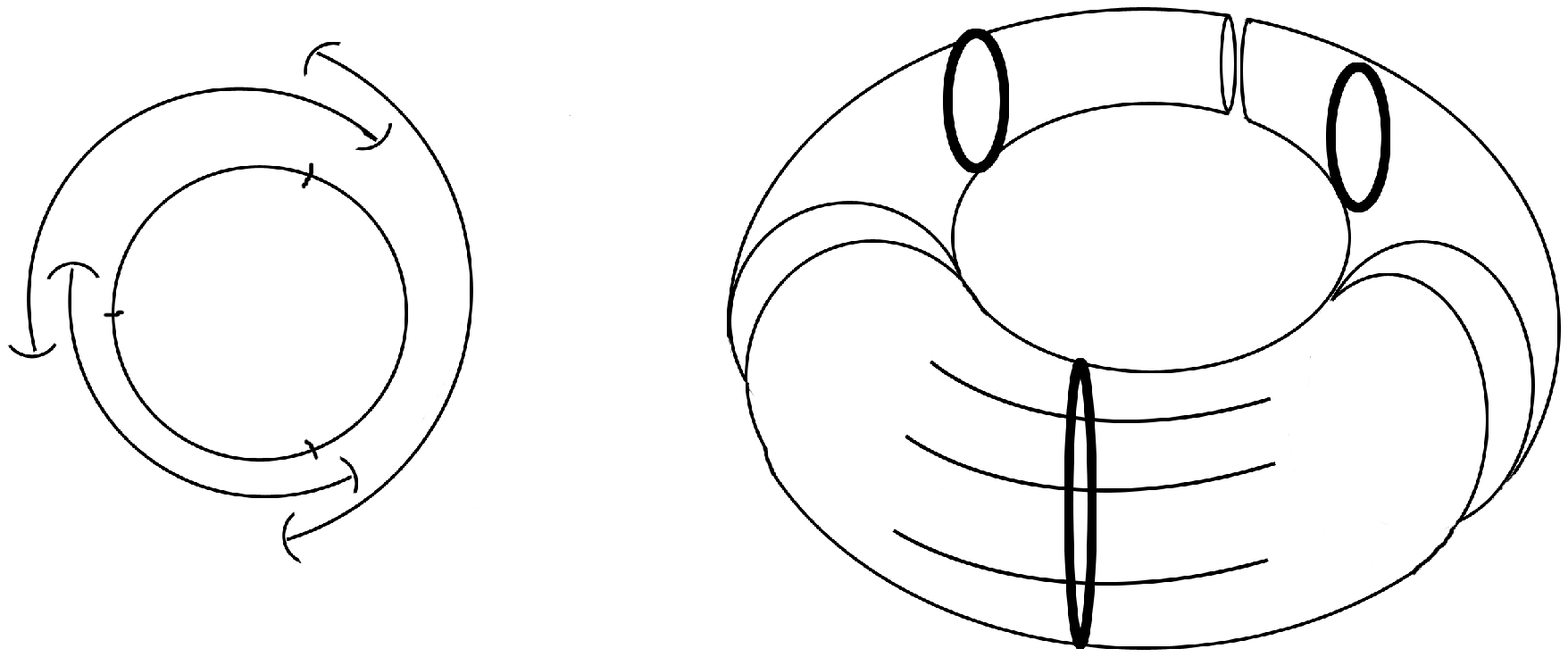}
\end{center}
\end{example}
One then takes a partition of unity $\{\rho_\alpha\}$ subordinate to $\mathcal{U}$ and writes $\tau \in \Omega^k(M)$ as $\Sigma \rho_\alpha \tau$. The requirement of $\mathcal{U}$ being a regular cover makes $\fol|_{U_\alpha}$ a fiber bundle with contractible fiber, and we integrate each $\rho_\alpha \tau$ by the standard integration over fibers. The result is a smooth form on $T$ of degree $k-dim(\fol)$. The main content of Theorem 3.1 of [Hae] is the following
\begin{theorem} Two choices of partition of unity give rise to forms on $T$ which are equal after pullback by holonomy. \qed
\end{theorem}
The kernel of the map $\int_\fol$ is then identified:
\begin{definition}
A $(p+k)$-form $\alpha \in \Omega^{p+k}(M)$ is said to be \textbf{ $\fol$-trivial} if
$$\alpha(X_1,\ldots X_p,Y_1,\ldots Y_k)=0\ \ \ \ \ \ \ \ \forall\ X_i \in T\fol,\ i =1,\ldots p.$$
\end{definition}
\begin{theorem} ([Hae], \S3.2)
The kernel of $\int_\fol$ is (the vector subspace generated by) $\fol$-trivial forms and differentials of $\fol$-trivial forms. 
\end{theorem}
From this, it follows that in the case $k=0$, the above map yields an isomorphism
$$H^{dim(\fol)}_\fol(M) \stackrel{\cong}{\longrightarrow} \Omega^0_c(T/H).$$

Consider our example $Y_\alpha$. The first obstruction $[\{\G,\G\}]$ lies in $H^{dim(\fol)}_\fol(M)$, and so we may apply the above apparatus. It frees us from the restriction of trying to solve $d_\fol \Sigma = \{\G,\G\}$ by integrating globally. Instead we integrate locally, and then ask if the solutions on different transverse patches "glue," in the above weaker sense. We now carry this out.

In the example, a regular cover is unnecessary in the $X_{H_3}$ direction, as the foliation is already a fiber bundle. In the $X_{H_\alpha}$ direction, on one solid torus, say $(r_1,\theta_1,\theta_2)$, we use the regular cover given in the example above. The resulting complete transversal is then three transverse disks, call them $T_1, T_2, T_3$. Any function $F$ on one transverse disk may be written as $G - h^*G$ as follows, for $h$ a diffeomorphism between two disks. Call $F|_{T_1} \equiv F_1$. Then we solve
$$F_1 = G_1 - h_{12}^*G_2$$
by simply setting $G_1 = F_1$, $G_2 \equiv 0$. In fact we may always solve for a composition of two pullbacks, and it is only when considering a triple pullback $T_1 \to T_3 \to T_2 \to T_1$ from $T_1$ to itself, representing non-trivial holonomy, that we find an obstruction. Thus Haefliger's map in this case can be understood as picking one transverse disk $\theta_2=0$, integrating along the foliation from $\theta_2=0$ to $\theta_2 = 2\pi$, and asking if the resulting function can be written in the form $G-h^*G$ for $h$ a holonomy self-diffeomorphism of the disk.

We proceed in our obstructed cases by picking $\Gamma$ such that $\{\Gamma,\Gamma\}$ can be written without reference to the foliation coordinates, so $\{f,g\}_P$ makes sense, and integrating over $\fol$ only serves to multiply by a function which does not depend on the angular transverse coordinate.

\section{Obstructedness}

In the case $\alpha \in \Q$, leaf topology will cause obstructedness via a zero-average condition. In the case $\alpha \notin \Q$, while leaf topology is lost, there is another significant restriction, that of the leaves being dense in tori. In particular, this implies the zero average condition is retained. In the Liouville case, this will turn out to be sufficient to retain obstructedness. In the complementary diophantine case, however, another phenomenon takes over: the holonomy also works to kill the transverse structure, making it impossible to support a non-trivial Poisson bracket. Thus we will see how the foliation may act to revive unobstructedness.

As we prove vanishing or non-vanishing of $[\{\G,\G\}]$ via Haefliger's isomorphism, we find that not only is transverse structure retained in the Liouville case in spite of significant holonomy restrictions, but in addition, Haefliger's map successfully transfers this structure to a fine model for the leaf space.

\subsection{Construction of $\Gamma$, $\alpha \in \mathbb{Q}$ case}

We produce $f,g \in C^\infty(Y_\alpha)$ such that $\{f,g\}_P$ is not in the denominator of $\Omega^0(T/H)$. For simplicity we refer to $r_1$ as $r$, and begin by taking $\rho(r)$ a smooth cutoff function vanishing smoothly on a neighborhood of both "poles" $r = 0, \frac{1}{\sqrt{2}}$, allowing us to consider only one transverse disk.  Furthermore choose $\rho(r)=r$ over some interval away from the poles. Suppose $\alpha = \frac{p}{q}$. Then consider
$$\begin{array}{rcl}
f(r,\theta_1) &:=& \rho(r)sin(2\pi q\theta_1),\\
g(r,\theta_1) &:=& \rho(r),\\
\end{array}$$
so on the corresponding annulus, 
$$\Gamma = r\cdot sin(2\pi q\theta_1)(X_{H_\alpha})^* + r(X_{H_3})^*.$$
We note that $\rho$ does not affect global closedness, as $d_\fol$ does not involve the radius at all, and $\Gamma$ so defined is easily checked to be $d_\fol$-closed. We compute $\{f,g\}_P$, where in terms of the transverse coordinates $(r,\theta_1)$, $\omega|_{Y_\alpha} = rdr \land d\theta_1$, so

$$\begin{array}{rcl}
\{f,g\}_P &=& \frac{1}{r}(\frac{\partial f}{\partial r}\frac{\partial g}{\partial \theta_1} - \frac{\partial f}{\partial \theta_1}\frac{\partial g}{\partial r})\\
&=& -2\pi qcos(q\theta_1)\\
\end{array}$$
Thus we must show that $F(r,\theta_1) = qcos(q\theta_1)$ is not of the form $G-h^*G = G(r,\theta_1) - G(r,\theta_1 + \alpha)$ for any smooth function $G$ on the disk. Suppose it were, and consider the equation in terms of Fourier series. That is, writing $F(r,\theta_1) = \Sigma_{\substack{n \in \Z}} F_n(r)e^{2\pi in\theta_1}$, etc., we have:

$$\begin{array}{rcl}
F(r,\theta_1) &=& G(r,\theta_1) - G(r,\theta_1 + \alpha)\\
\frac{q}{2}e^{2\pi i q\theta_1} + \frac{q}{2}e^{-2\pi i q\theta_1} &=& \Sigma G_n(r) e^{2\pi in\theta_1} - G_n(r) e^{2\pi in(\theta_1 + \alpha)}\\
&=& \Sigma G_n(r) e^{2\pi in\theta_1}(1- e^{2\pi in\alpha})\\
\end{array}$$
so we obtain 

$$\frac{q}{2} = 0,$$
as $q\alpha \in \mathbb{Z}$. Thus there is no solution, even for a single value of $r$, and $Y_\alpha$ in the case $\alpha \in \Q$ is obstructed.\qed

\subsection{Construction of $\Gamma,\ \alpha$ Liouville case}

The proof will be an extension of the idea in the rational case. Note that for $\alpha$ any irrational number, we can solve the functional equation $F = G - h^*G$ by Fourier series:

$$\begin{array}{rcl}
F(r,\theta) &=& G(r,\theta) - G(r, \theta + \alpha)\\
F_n(r)e^{2\pi i n \theta} &=& G_n(r)e^{2\pi i n \theta} - G_n(r)e^{2\pi i n (\theta+\alpha)}\\
&=& G_n(r)e^{2\pi i n \theta}(1 - e^{2\pi i n \alpha})\\
\end{array}.$$
Hence
$$G_n(r) = \frac{F_n(r)}{1 - e^{2\pi i n \alpha}}\ \forall\ n \neq 0,$$
which defines a solution for any fixed $r$ away from zero if $\alpha \notin \mathbb{Q}$, as long as $F_0(r) := \int_{S^1(r)}F(r,\theta)d\theta = 0$.

What remains in question is if these new Fourier series $\Sigma G_ne^{2\pi in\theta_1}$ define a smooth function, and for certain $\alpha$, the so-called Liouville numbers, there exist smooth functions $F$ for which they do not. This comes directly from replacing the term $1-e^{2\pi iq\alpha}=0$ in the rational case with a sequence of very small numbers, and when we divide to solve for $G_n(r)$, we encounter the so-called small divisor problem causing the Fourier coefficients to blow up. All application of this phenomenon in the present work is standard - only its context is novel.

\subsubsection{Definition of Liouville numbers}

Suppose we have a sequence of rational numbers $\{\frac{p_i}{q_i}\}$ converging to $\alpha$, for instance the $n^{th}$ decimal approximation to $\alpha$. We may say this example converges at the ``rate'' of $\frac{1}{q_n}$, as $\frac{1}{q_n} = \frac{1}{10^n}$, and the remainder term can certainly be no larger than this. It may be shown that any irrational number may in fact be approximated at the rate of $\frac{1}{q_n^2}$. If we ask for convergence at the rate of $\frac{1}{q_n^k}$ for all $k$, however, we find that not every number has such an approximation, and the set of numbers which do is of Lebesgue measure zero - these are called the Liouville numbers.

\begin{definition}
Let $\alpha \in \mathbb{R}, \alpha \notin \mathbb{Q}$. We call $\alpha\ \mathbf{Liouville}$ if it has an approximation $\{\frac{p_n}{q_n}\}$ such that, for any $k \geq 1$, there exists a constant $\lambda(k)$ such that 
$$|\frac{p_n}{q_n} - \alpha| < \frac{\lambda}{q_n^k}$$
for all $n$, or equivalently that
$$|p_n - q_n\alpha| < \frac{\lambda}{q_n^k}$$
for all $n$.
\end{definition}

It is also equivalent to say $\alpha$ is Liouville if $\forall k\ \exists (p,q)$ such that
$$|\frac{p}{q} - \alpha| < \frac{1}{|q|^k}.$$
Specifically it is not directly a question of convergence in $n$.
\begin{example}
Liouville's original example was the following: $\sum_{\substack{i=1}}^{\substack{\infty}} \frac{1}{10^{n!}}$.
\end{example}
Recall the strikingly similar fact that $F = \sum F_ne^{2\pi in\theta}$ is in $C^\infty(S^1)$ iff for any $k \geq 1$, there exists a constant $\lambda(k)$ such that 
$$|F_j| < \frac{\lambda}{j^k} \text{ for all } j.$$
Given a sequence $\{\frac{p_n}{q_n}\}$, if we consider a Fourier series having non-zero coefficients only at $j=q_n$, the smoothness condition becomes
$$|F_{q_n}| < \frac{\lambda}{q_n^k} \text{ for all } n.$$

Define $F$ by setting its Fourier coefficients $F_{q_n} := p_n - q_n \alpha$, and zero otherwise. Then $F$ is smooth, due to the Liouville condition for $\alpha$. However, using the geometric fact that a chord on a circle is shorter than the arc subtended, in the form of the inequality
$$|1-e^{2\pi i(p-q\alpha)}| = |1-e^{2\pi iq\alpha}| \leq 2\pi|p - q\alpha|,$$
we see that the coefficients $G_j$ satisfy the following:
$$|G_n| = \frac{|F_n|}{|1-e^{2\pi i n\theta}|} \geq \frac{|F_n|}{2\pi|p_n - q_n\alpha|} = \frac{1}{2\pi},$$
and so do not even give a series whose terms converge to zero.

We obtain a function on the disk by multiplying $F(\theta)$ by the cutoff function $\rho(r)$ used in the rational case. Consider $f(r,\theta) := \rho(r)F(\theta), g(r,\theta) := \rho(r)r$. We compute, restricting our attention to the annulus where $\rho(r)$ is identically $r$, and find

$$\begin{array}{rcl}
\{f,g\}_P &=& \{r\Sigma\ F_nsin(2\pi n\theta), r\}_P\\ \\
&=& \Sigma\ F_n\{rsin(2\pi n\theta), r\}_P\\ \\
&=& -2\pi\Sigma\ nF_ncos(2\pi n\theta)\\
\end{array}$$
This Fourier series has the same property as $F$, namely that it converges to a smooth function, yet by the above argument, its corresponding $G$ does not, even for one value of $r$. Therefore $\li_2(\Gamma, \Gamma)$ is not sent to zero in $\Omega^0(T/ H)$ under $\int_\fol$, and this infinitesimal deformation is obstructed. \qed

\section{Unobstructedness}

We may understand the obstructedness result in the Liouville case above as a result of the zero-average condition, brought on by dense leaves intertwining the tangent and transverse structure. We now proceed to show how such holonomy may produce a contrary effect, killing the transverse topology necessary to support an obstruction class, and thereby reviving unobstructedness.

\begin{definition}
A number $\alpha \in \mathbb{R}$ is said to satisfy a $\mathbf{diophantine\ condition}$ if it is not Liouville, that is, if there exists a $k \geq 0$ such that for any non-zero integers $p$ and $q$,
$$|\frac{p}{q} - \alpha| > \frac{1}{q^k}$$.
\end{definition}

We prove that in the case $\alpha$ is diophantine, every infinitesimal deformation is the infinitesimal deformation of some full deformation (i.e. the dimension of the moduli space is the same as the dimension of its tangent space). The final result is formal, but may integrate to a $C^\infty$ result as well. 

To establish formal unobstructedness, the first obstruction we must prove to vanish is the Kuranishi map. This is where we focus our attention, and in the final section we show how higher obstructions reduce to this case.

Let
$$\Gamma = fX^*_{H_\alpha} + gX^*_{H_3} \in \Omega^1(\fol),$$
such that $d_\fol \Gamma = 0$. Then we seek to show there exists
$$\Sigma = h_1X^*_{H_\alpha} + h_2X^*_{H_3} \in \Omega^1(\fol)$$ such that
$$d_\fol \Sigma = \{\Gamma, \Gamma\}.$$
(Recall that $f,g,$ and $h_i$ are all functions in $C^\infty(Y_\alpha)$.) In earlier discussion we saw that this amounts to showing $\exists\ h_1, h_2 \in C^\infty(Y_\alpha)$ such that
$$X_{H_\alpha}(h_2) - X_{H_3}(h_1) = \{f,g\}_P.$$
In fact we solve
$$X_{H_\alpha}(h) = \{ \tilde{f},\tilde{g}\}_P,$$
with modified versions $\tilde{f},\tilde{g}$ of $f,g$. Recall that the Kuranishi map is well-defined on cohomology, so we may first modify $\Gamma$ by an exact 1-form $d_\fol h$.

\subsection{Main proposition}

Such modification by an exact form will allow us to assume the function on the right-hand side has a special property: $\{\tilde{f},\tilde{g}\}_{0,0} = 0$ $\forall \ r \neq 0$. Here $\{\tilde{f},\tilde{g}\}_{0,0}$ is the constant term in the double Fourier series for $\{\tilde{f},\tilde{g}\}$ in $\theta_1, \theta_2$, and we have the following, our main proposition: 

\begin{proposition}\label{mainProp}Given $\phi:Y_\alpha \to \R$ such that $\phi_{0,0} \equiv 0 \ \forall \ r \neq 0$, the DE
$$X_{H_\alpha}(h) = \phi$$
always has a smooth solution $h$ on the open set $r \in [0,\frac{1}{\sqrt{2\alpha}})$. Furthermore, this solution is unique up to an even function of one variable.
\end{proposition}
Equivalently, we may always solve
$$X_{H_\alpha}(h) = -(\phi - \phi_{0,0}).$$
We assume this proposition for the moment, and proceed to prove our result for the case of flat transverse connection $G$. We return later to address the case of non-flat $G$. Recall that $f_{0,0}$, and thus also $\frac{\partial f_{0,0}}{\partial r}$ as well, are constant in $\theta_1$ and $\theta_2$. Then denoting $\tilde{g} \equiv g + X_{H_3}(h)$, we have
$$\begin{array}{rcl}
\{\Gamma + d_\fol h, \Gamma + d_\fol h\} &=& \{f+X_{H_\alpha}(h),\ \tilde{g} \}\ ((X_{H_\alpha})^* \land (X_{H_3})^*)\\
\\
\{f+X_{H_\alpha}(h),\ \tilde{g} \}_{0,0} &=& \{f_{0,0},\ \tilde{g}\}_{0,0}\\
&=& \int \int \frac{1}{r} (\frac{\partial f_{0,0}}{\partial r}\ \frac{\partial\tilde{g}}{\partial \theta_1} - \frac{\partial f_{0,0}}{\partial \theta_1}\frac{\partial\tilde{g}}{\partial r}) d\theta_1d\theta_2\\
&=& \frac{1}{r}\frac{\partial f_{0,0}}{\partial r} \int \int \frac{\partial\tilde{g}}{\partial \theta_1}\ d\theta_1d\theta_2\\
&=& 0.
\end{array}$$

But now we simply need to solve the same DE again, as it suffices to solve
$$X_{H_\alpha}(h) = \{ \tilde{f},\tilde{g}\}_P$$
where by the above, we know that $\{ \tilde{f},\tilde{g}\}_{0,0}(r) = 0\ \forall \ r \neq 0$.
\qed

\subsubsection{Retrospect}

This result may be viewed in light of basic functions $\Omega^0_{bas}$ and Haefliger's group $\Omega^0_c(T/H)$. In the irrational case in general, $\Omega^0_{bas}$ is some subset of smooth functions of one variable, $r$, for basic functions are constant on tori. In the diophantine case in particular, the above says we can always solve the exactness equation given zero average, so $H^{top}_\fol \cong \Omega^0(T/ H) \cong \Omega^0_{bas}$ (the average is a basic function). Thus we are essentially carrying out all of our computations on a one-dimensional space, which cannot support a Poisson bracket. In the Liouville case, $\Omega^0_c(T/H)$ is in fact infinite-dimensional as a module over $\Omega^0_{bas}$, and plays the role of a de facto second variable.

\subsection{Proof of proposition}

\subsubsection{Solution away from r=0}

Consider the differential equation
$$\frac{\partial h}{\partial \theta_1} + \alpha \frac{\partial h}{\partial \theta_2} = f$$
for $h,f: Y_\alpha \to \R.$ As long as $r_1 \in (0, \frac{1}{\sqrt{2}})$, this may be written in terms of Fourier series in the two angular coordinates $\theta_1, \theta_2$. We mention explicitly that all of the coefficients are functions in $(r, \theta_3)$. Then using the symbols $\xi_n = e^{2\pi i \theta_n}$, $n=1,2$, we have
$$\begin{array}{rcl}
f_{p,q}\ \xi_1^p\xi_2^q &=& X_{H_\alpha}(h_{p,q}\ \xi_1^p\xi_2^q)\\
&=& 2\pi i (p + \alpha q) h_{p,q}\ \xi_1^p\xi_2^q.
\end{array}$$
Smoothness of $f$ implies rapid decay of the Fourier coefficients: $\forall\ k \geq 1, \exists\ \lambda (k)$ such that
$$|f_{p,q}| < \frac{\lambda}{|q|^k}\ \ \forall\ q.$$
By the diophantine condition on $\alpha$, $\exists\ j$ such that $\forall\ (p,q)$, 
$$|p + q \alpha| > \frac{1}{|q|^j}.$$
Combining these, we obtain that $\forall\ k \geq j +1, \exists\ \lambda (k)$ such that
$$|h_{p,q}| = \frac{|f_{p,q}|}{ 2\pi\ |p + \alpha q|} < \frac{\lambda }{ 2\pi\ |q|^k}\ \ \forall\ q.$$
The inequality follows for all $k$ by modification of a finite number of constants. Thus $h_{p,q}$ necessarily defines the Fourier coefficients of a smooth function, given only the obvious precondition that $f_{0,0}=0$. Now we have solutions on each torus, and since the coefficients $f_{p,q}$ vary smoothly in $r$, the solutions vary smoothly in $r$, and we have a smooth solution everywhere but r=0.

\subsubsection{Alternative view of Haefliger's isomorphism - Examples}

In the obstructed cases, we chose forms $\G$ which could be written locally without reference to the foliation variables, in order to make use of the isomorphism $H_\fol^{dim\fol} \cong \Omega^0(T/H)$. Here we are given arbitrary forms, but wish to make use of this tool all the same. 

We first distinguish the two issues complicating solution of the PDE in the proposition: the ``singularity'' at $r=0$, and the Liouville phenomenon on tori.
\begin{example}\label{example1} Consider the PDE on the disk $D^2$
$$\frac{\partial h}{\partial \theta} = f(r,\theta),$$
where $f$ is an arbitrary given smooth function on $D^2$. We consider the smooth blowup of $D^2$ at the origin (Ch.1, \S2 of [Ar]), viewed as the Mobius strip with coordinates $(x,\theta)$ in the rectangle $(-1,1) \times [0,\pi]$, identifying $(x,0) \sim (-x,\pi),$ and mapping to the disk $(r,\theta_1)$ by
\begin{equation}
\label{eq:mdiv}
(r,\theta_1) \impliedby
\begin{cases}
(x,\theta) & \text{if $x > 0$} \\
(|x|, \theta + \pi) & \text{if $x < 0$}
\end{cases}
\end{equation}
\begin{center}
\includegraphics{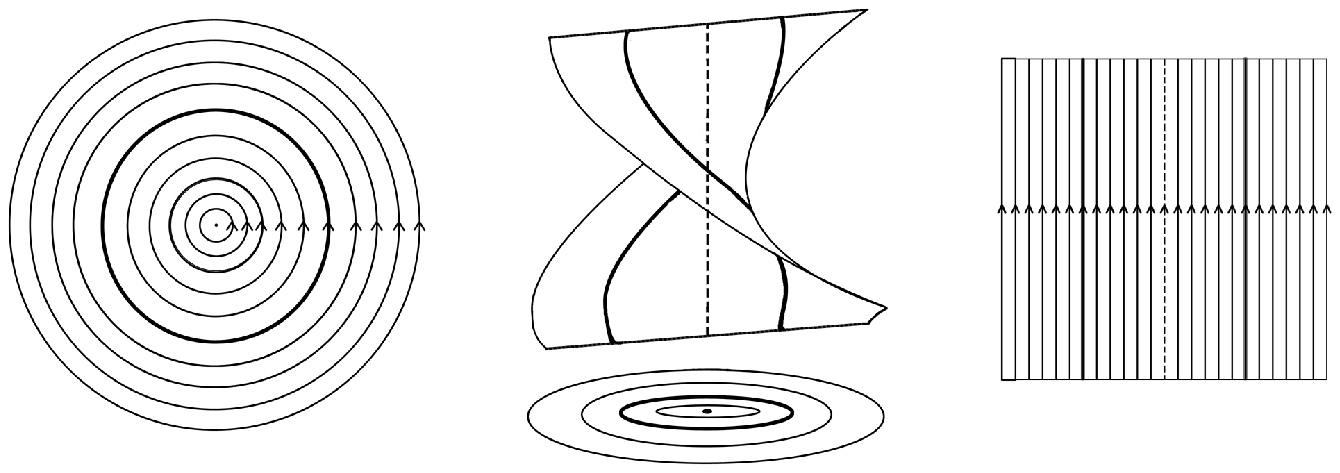}
\end{center}
On the blowup, the PDE is non-degenerate, and we denote the lift by
$$\frac{\partial \tilde{h}}{\partial \theta} = \tilde{f}(x,\theta),$$
where $f(0)$ lifts as constant in $\theta$. This may clearly be solved on the rectangle, uniquely up to a function $C(x)$,
$$\tilde{h}(x,\theta) := \int_0^\theta \tilde{f} d\theta + C(x).$$
Passing to the Mobius strip, this is further restricted by the matching condition
$$\begin{array}{rcl}
\tilde{h}(x,0) &=& \tilde{h}(-x,\pi)\\
C(x) &=& \int_0^\pi \tilde{f}(-x)d\theta + C(-x).\\
\end{array}$$
So $C(x)$ is in fact partially restricted, with its odd part determined by
$$\frac{C(x)-C(-x)}{2} = \frac{1}{2}\int_0^\pi \tilde{f}(-x)d\theta.$$
To achieve a sufficient condition, we see already on the disk that solving on any circle $S^1(r)$ has the necessary and sufficient condition
$$\int_{S^1(r)} f = 0.$$
This condition lifts and extends continuously to $x=0$, where it implies that $\tilde{f}(0,\theta) \equiv 0,$ and thus that $\tilde{h}$ is constant along $x=0$, and so descends back to $D^2$.
\end{example}

This calculation may be interpreted as a solution of $\Gamma = d_\fol h$ on the disk foliated by circles, where $\Gamma = f(\frac{\partial}{\partial \theta})^*.$ Thus we have determined that a form is $d_\fol$-exact iff it satisfies the zero average condition $\forall r \neq 0$. This may be interpreted instead as oddness of the function $$g(x) := \int_0^\pi \tilde{f}(x, \theta)d\theta,$$ and so deduce that
$$H_\fol^1(D^2) \cong \{even\ functions\ on\ (-1,1)\},$$
which we may understand as smooth functions on the orbifold $(-1,1)/x\sim-x$.
\\

We now revisit the above discussion of Fourier series on the torus using this perspective.

\begin{example}\label{example2} 
Consider the DE on the torus $(\theta_1,\theta_2)$:
$$\begin{array}{rcl}X_{H_\alpha}(h) &=& f\\
(\frac{\partial}{\partial \theta_1} + \alpha \frac{\partial}{\partial \theta_2})(h) &=& f(\theta_1,\theta_2).\\
\end{array}$$
By the standard method of characteristics for linear, first-order PDEs (Ch.2, \S7 of [Ar]), this is equivalent to a system of total differentials:
$$d\theta_1 = \frac{d\theta_2}{\alpha} = \frac{dh}{-f(\theta_1,\theta_2)}$$
(More precisely, each is the pullback of the named differential to $\R$ along the inclusion of any orbit of $X_{H_\alpha}$.) We integrate along each characteristic curve, yielding from the first equation the defining equation for the foliation:
$$\theta_2 = \alpha \theta_1 + const.$$
Using the second equality we solve for $h(\theta_1,\theta_2)$:
$$h(\theta_1,\theta_2) := \frac{1}{\alpha}\int_0^{\theta_2} f(\theta_1,s) ds + C(\theta_1),$$
up to the function $C(\theta_1)$. As in the previous example, $C(\theta_1)$ is not arbitrary, but is restricted by the matching, or holonomy condition:
$$\begin{array}{rcl}
h(\theta_1,0) &=& h(\theta_1 + \alpha2\pi, 2\pi)\\
C(\theta_1) &=& \frac{1}{\alpha}\int_0^{2\pi} f(\theta_1,s) ds + C(\theta_1 + \alpha2\pi)\\
\end{array}$$
Denoting $F(\theta_1) := \frac{1}{\alpha}\int_0^{2\pi} f(\theta_1,s) ds$, we note that any smooth function $F(\theta_1)$ on the circle can be thus obtained, e.g. from $f(\theta_1,\theta_2) = \frac{1}{2\pi}F(\theta_1)$, and so this is exactly the homological equation of a deformation of a rotation,
$$F(\theta_1) = C(\theta_1)-C(\theta_1 + \alpha2\pi).$$
In contrast to decomposition into even and odd functions, a decomposition of an arbitrary function by this holonomy condition is not always possible, as seen in the Liouville case. However, in the diophantine case, the necessary condition $\int_{S^1}F(\theta_1)d\theta_1 = 0$ is sufficient. Interpreting this as exactness in foliation cohomology for the linear foliation $\fol_\alpha$ of $T^2$ of slope $\alpha$, $H^1_\fol(T^2, \fol_\alpha)$ is parametrized by
$$\int F d\theta_1 = \int \int f d\theta_1 d\theta_2 = f_{0,0},$$
and for $\alpha$ diophantine, we recover the well-known result about circle diffeomorphisms (see, e.g., [dS]):
$$H^1_\fol(T^2, \fol_\alpha) \cong \R.$$
\end{example}

\subsubsection{Solid torus - local proof of proposition}

We return to the proof of Proposition \ref{mainProp}, combining the previous two examples. The differential equation of the proposition is the same as in Example \ref{example2}, only depending nominally on $r$ as well:
$$X_{H_\alpha}(h) = f(r,\theta_1,\theta_2).$$
Integrating over $\theta_2$ as above now puts us on a transverse disk $(r,\theta_1)$. The homological equation is likewise the same, extended to the disk:
$$F(r,\theta_1) = C(r,\theta_1)-C(r,\theta_1 + \alpha2\pi).$$
First we solve again away from $r=0$, only now in the context of the homological equation. We saw in the earlier analysis of Fourier series that $f_{0,0}(r) = 0$ is a necessary and sufficient condition for solution of the DE for any $r\neq0$, and in Example \ref{example2} above we saw how integrating over $\theta_2$ translates that condition into zero average of $F(r,\theta_1)$ over each circle $S^1(r)$. The Fourier series version of this homological equation was already examined as well, in the Liouville case:
$$C_n(r) = \frac{F_n(r)}{1 - e^{2\pi i n \alpha}}.$$
Here we use the opposite comparison to that used in the Liouville case,
$$|1-e^{2\pi in\alpha}| \geq \pi|p_n - n\alpha|$$
for some $p_n$, namely the closest integer to $n\alpha$. Combining the decay of $|F_n|$ and the diophantine condition on $\alpha,$ we obtain that $\exists j$ such that $\forall\ k \geq j+1, \exists\ \lambda (k)$
$$|C_n(r)| = \frac{|F_n(r)|}{|1 - e^{2\pi i n \alpha}|} \leq \frac{|F_n(r)|}{\pi|p_n - n\alpha|} < \frac{\lambda}{\pi |n|^k}.$$
Again, we obtain the inequality for all $k$ by modification of a finite number of constants.

To solve the homological equation on the disk, we emulate the argument of Example \ref{example1}, only now applied to an algebraic question. We first lift the homological equation, not to the blowup, but to its double cover. Namely, we write
$$\tilde{F}(x,\theta_1) = \tilde{C}(x,\theta_1)-\tilde{C}(x, \theta_1 + \alpha2\pi),$$
where $\tilde{F}$ is canonically extended to the double cover $(-1,1) \times [0,2\pi]$, with $(x,0) \sim (x,2\pi)$. The necessary condition of zero average over circles again lifts, and extends continuously to $x=0$. Since $\tilde{F}(0,\theta)$ is in fact constant, zero average implies that constant is zero. Given zero average for all $x$, a solution to the homological equation clearly extends smoothly across $x=0$. In fact, since $\tilde{F}(0,\theta_1) \equiv 0$,
$$\tilde{C}(0,\theta_1)= \tilde{C}(0, \theta_1 + \alpha2\pi).$$
But now, the orbit of a point under repeated translation by $2\pi \alpha$ is dense, and $\tilde{C}(0,\theta_1)$ must be constant on such an orbit, so it must be constant, and the solution descends, not just to the Mobius band, but in fact to the disk. Thus the proposition is proven, for either solid torus, with a flat transverse connection.
\qed
\\

We note that the equation $\Gamma = d_\fol h$ can be solved under the same conditions as in the plain disk case, namely for zero average over circles in the homological equation, which corresponds in the original differential problem to zero average over tori. In this case, we have the representation of $H^2_\fol$ as functions on the orbifold $(-1,1)/x \sim -x$ as follows. To pass to the Mobius strip from its double cover, the Fourier coefficients of $\tilde{F}$ must satisfy
$$\begin{array}{rcl}
\tilde{F}(x,\theta_1) &=& \tilde{F}(-x,\theta_1+\pi)\\
\tilde{F}_n(x)e^{in\theta_1} &=& \tilde{F}_n(-x)e^{2\pi in\theta_1}\cdot e^{2\pi in\pi}\\
&=& \tilde{F}_n(-x)e^{2\pi in\theta_1}\cdot (-1)^n\\
\end{array}$$
So as a function of $x$, $\tilde{F}$ is arbitrary for $x \in (0,1)$, $\tilde{F}(x)$ determines $\tilde{F}(-x)$, and at $x=0$, the $\tilde{F}_n$ are the coefficients of an arbitrary even function.

\subsubsection{Non-zero curvature - global proof}
Now that we have proven exactness of the first obstruction for a flat connection on each solid torus, we extend the result to all of $Y_\alpha$, sacrificing flatness. We specify a connection by picking a distribution $G$ in $TY_\alpha$ complementary to $T\fol$. The distribution used for the flat connection on each solid torus $(r_i, \theta_1, \theta_2)$ was $G_i := span\{\frac{\partial}{\partial r_i}, \frac{\partial}{\partial \theta_i}\}$, for $i=1,2$. We define a global $G$ by considering an $\epsilon$-neighborhood of $r_1 = \frac{1}{2\sqrt{2}}$, setting $G=G_i$ to either side of that neighborhood, and on the neighborhood make use of a smooth homotopy from $\frac{\partial}{\partial \theta_1}$ to $\frac{\partial}{\partial \theta_2}$. The span of the other basis vector $\frac{\partial}{\partial r_i}$ is the same across both patches.

Recall the problem to prove unobstructedness was to solve for $\Sigma$, or alternatively $h$:
$$\begin{array}{rcl}
d_\fol \Sigma &=& \{\Gamma, \Gamma\}\\
X_{H_\alpha}(h) &=& \{ \tilde{f},\tilde{g}\}_P\\
 &=& P^{ij}(\nabla_i(\tilde{f})\nabla_j(\tilde{g}) - \nabla_j(\tilde{f})\nabla_i(\tilde{g}))\\
\end{array}$$
Recalling the expression for $\nabla$ in coordinates,
$$\nabla_i\Gamma_\alpha = \frac{\partial\G_\alpha}{\partial y^i} + R^\beta_i\frac{\partial \G_\alpha}{\partial q^\beta} + \G_\beta \frac{\partial R^\beta_i}{\partial q^\alpha}$$
it is evidently desirable to choose $G$ such that $R_i^\beta$ does not depend on the foliation coordinates $q^\alpha$. For $t \in [0, \frac{1}{\sqrt{2}}]$, let $\tau(t)$ be a smooth function such that
\begin{equation}
\tau(t)=
\begin{cases}
    0 & \text{for $t \leq \frac{1}{2\sqrt{2}}-\epsilon$,}\\
    1 & \text{for $t \geq \frac{1}{2\sqrt{2}}+\epsilon$,}\\
\end{cases}
\end{equation}
define
$$X(r_1) := (1-\tau(r_1))\frac{\partial}{\partial \theta_1} - \alpha\tau(r_1)\frac{\partial}{\partial \theta_2},$$
and set $G = span\{\frac{\partial}{\partial r_1}, X(r_1)\}$. 
Polar coordinates on a transverse disk are still natural \textit{local} transverse coordinates:
$$\begin{array}{rcl}
y^1 = r_1&,& \frac{\partial}{\partial q^1} = X_{H_\alpha},\\
y^2 = \theta_1&,& \frac{\partial}{\partial q^2} = \frac{\partial}{\partial \theta_3}\\
\end{array}$$
In terms of these our basis element $e_2 = X$ for $G$ is given in one chart by:
$$X = \frac{\partial}{\partial \theta_1} - \tau(r_1)X_{H_\alpha}$$
so $R^1_2 = -\tau(r_1)$, with all other coefficients zero. Similarly, denoting by $\tilde{\tau}$ a suitable reparametrization of $\tau$ to use $r_2$ as the variable, we have in the other chart
$$X = \alpha \frac{\partial}{\partial \theta_2} + (1 - \tilde{\tau}(r_2))X_{H_\alpha}.$$

We wish to use Proposition \ref{mainProp}, and are thus trying to show $\{f_{0,0},\tilde{g}\}_{0,0} = 0$. In other words,
$$\begin{array}{rcl}
0 &=& \int \int P^{ij}(\nabla_i(f_{0,0})\nabla_j(\tilde{g}) - \nabla_j(f_{0,0})\nabla_i(\tilde{g})) d\theta_1 d\theta_2\\
&=& \frac{1}{r}\int \int\ \frac{\partial f_{0,0}}{\partial r_1}(\frac{\partial} {\partial \theta_1} - \tau(r_1) X_{H_\alpha})(\tilde{g})\\
& &\ \ \ \ \ \ \ \ \ - (\frac{\partial} {\partial \theta_1} - \tau(r_1) X_{H_\alpha})(f_{0,0})\frac{\partial \tilde{g}}{\partial r_1}\ d\theta_1 d\theta_2\\
\end{array}$$
The first and third terms vanish as before. The fourth term vanishes as $X_{H_\alpha}(f_{0,0}) = 0$. Thus we are left with

$$\int \int -\frac{\partial f_{0,0}}{\partial r_1} \tau(r_1) X_\alpha(\tilde{g})  d\theta_1 d\theta_2 = -\frac{\partial f_{0,0}}{\partial r_1} \tau(r_1) \int \int  \frac{\partial \tilde{g}}{\partial \theta_1} + \alpha  \frac{\partial \tilde{g}}{\partial \theta_2} d\theta_1 d\theta_2 = 0$$
\\

We have shown we can solve locally even for this non-flat $G$. As the solution is unique after choosing an even function, we solve on each solid torus, choosing the zero function. Then since $G$ is global, if we consider two open solid tori $r_1 \in [0,\frac{1}{2\sqrt{2}}+\epsilon)$ and $r_1 \in (\frac{1}{2\sqrt{2}}-\epsilon, \frac{1}{\sqrt{2}}]$, the solutions must agree on the overlap, and thus patch to a global solution. As an existence proof, we have determined that there exists a solution on the quotient space $S^2/\sim$ iff the average over each circle $r_1 = const$ is zero, and the solution is unique up to a "holonomy-even" function on the whole space.

\subsection{Vanishing of higher obstructions}

We have as yet not given the definition of the higher $L_\infty$ operators, as they appear in [OP]. We recall the earlier claim that the higher-order operators of the $L_\infty$-structure vanish exactly when there exists a complementary foliation everywhere transverse to the null foliation. In coordinates we had
$$F_G = F_{ij}^\alpha \frac{\partial}{\partial q^\alpha} \otimes \ e_i^* \land e_j^*.$$
We denote by $F^\sharp := F_G\omega^{-1}$ the raising of indices, so $(F^\sharp)_i^{\alpha j} = F_i^{\alpha j} := F_{ik}^\alpha P^{kj}$, and $F^\sharp \in \G(End(N^*\fol)\otimes T\fol)$. For the global $G$ chosen above in the example $Y_\alpha$, this curvature takes a particularly simple form. The coordinate expression for the curvature is 
$$\begin{array}{rcl}
F^\alpha_{ij} &=& \frac{\partial R^\alpha_j}{\partial y^i} - \frac{\partial R^\alpha_i}{\partial y^j} + R^\beta_i \frac{\partial R^\alpha_j}{\partial q^\beta} - R^\beta_j \frac{\partial R^\alpha_i}{\partial q^\beta}\\
F^1_{12} &=& -\frac{\partial \tau}{\partial r_1}.
\end{array}$$
So denoting $J = \left( \begin{array}{cc}
0 & -1 \\
1 & 0 \end{array} \right)$, we have
$$\begin{array}{rcl}
(F^1_{ij}) &=&  \frac{\partial \tau}{\partial r_1} J,\\
(F^2_{ij}) &=& 0,\\ 
(F_i^{1j}) &=& (F^1_{ik}) P^{kj} = \frac{\partial \tau}{\partial r_1} J \frac{1}{r}(-J) = \frac{1}{r}\frac{\partial \tau}{\partial r_1} Id.
\end{array}$$

We are now ready to give the definition of the higher $L_\infty$ operators.
\begin{definition}
We define $\li_k: \Gamma(\tensor^kT^*\fol) \to \Gamma(\Lambda^2T^*\fol)$ for $k \geq 3$:
$$\li_k(\G_1,\ldots \G_k) := \Sigma_{\sigma \in S_k}(-1)^{|\sigma|}P(\nabla\G_{\sigma(1)}, (F^\sharp\rfloor\G_{\sigma(2)})\ldots(F^\sharp\rfloor\G_{\sigma(k-1)})\nabla\G_{\sigma(k)}).$$
For $\G_1 = ... = \G_k = \G$, this reduces to
$$\li_k(\G,\ldots \G) := P(\nabla\G, (F^\sharp\rfloor\G)^{k-2}\nabla\G).$$
\end{definition}

\vspace{.2in}
We make a few remarks pertinent to our example. The contraction $F^\sharp\rfloor\G$ ignores any $(\frac{\partial}{\partial \theta_3})^*$ component of $\G$. The resulting endomorphism of $N^*\fol$ is therefore
$$(F^\sharp\rfloor\G)^{k-2} = (\frac{1}{r}\frac{\partial \tau}{\partial r_1}f)^{k-2} Id.$$
Recalling further that we may assume $f = f_{0,0}$, we see that this entire expression pulls out of the double integral in the proposition, and the earlier arguments apply to the higher obstructions as well.
\vspace{.15in}

It only remains confirm that we may solve all of these individual exactness problems simultaneously and consistently. We first point out that when $G$ is integrable ($F\equiv 0$), all higher $L_\infty$ operators vanish, yet higher obstructions do still persist, namely the following:
$$\begin{array}{rcl}
t &:& d_\fol \Gamma_1 = 0;\\
\\
t^2 &:& -d_\fol \Gamma_2 = \frac{1}{2}\{\Gamma_1,\Gamma_1\};\\
\\
t^3 &:& -d_\fol \Gamma_3 = \frac{1}{2}(\{\Gamma_1,\Gamma_2\} + \{\Gamma_2,\Gamma_1\});\\
\\
t^4 &:& -d_\fol \Gamma_4 = \frac{1}{2}(\{\Gamma_1,\Gamma_3\} + \{\Gamma_3,\Gamma_1\} + \{\Gamma_2,\Gamma_2\});\\
\\
& &\ \ \  \vdots
\end{array}$$

We note that the main proposition does not depend on $f$ and $g$ coming from the same $\G$. For $\G_i := f_i (X_{H_\alpha})^* + g_i (X_{H_3})^*$,
$$\{\G_1,\G_2\} = P^{ij}(\nabla_i(f_1)\nabla_j(g_2) - \nabla_j(f_2)\nabla_i(g_1))(X_{H_\alpha}^* \land X_{H_3}^*).$$
We may modify $\G_1$ and $\G_2$ separately to make both $f_1$ and $f_2$ constant in $\theta_1$ and $\theta_2$, the crux of the argument. We modify all $\G_i$ in this manner. By the above argument, any higher bracket, all of whose arguments are $\G_i = f_i (X_{H_\alpha})^* + g_i (\frac{\partial}{\partial \theta_3})$ with $f_i$ constant in $\theta_1$ and $\theta_2$, is exact.

\qed

\section{Reflections on the minimal model}

In the above computations, we reduced a differential question about a foliation to an algebraic question on a complete transversal, as an \'etale model for the leaf space. The connection is rather tight - for instance the resolution of a singularity proved to be almost identical under the two perspectives.

The transverse bracket $\li_2$ in particular is well-represented by its image under Haefliger's map, which, due to the isomorphism $\Omega^0_c(T/H) \cong H^2_\fol(Y_\alpha)$ in this case, can be viewed as the induced bracket on $\li_1$-cohomology.

The transverse curvature used in the definition of the higher-order $L_\infty$-operators, on the other hand, was defined as a global vector-valued form, and in this example, it was simply shown to have no effect. A natural next step would be to consider an example with vanishing first obstruction, and use the curvature to twist the bracket such that higher obstructions do not vanish. Then the following natural questions arise: Does the curvature itself pass to an \'etale model for the leaf space, something of the form $\Omega^*(T/H;T\fol)$? Do the higher order $L_\infty$-operators descend? And what measures the difference, if any, between this $L_\infty$-structure and the induced dgLa structure on $\li_1$-cohomology?



\end{document}